\documentclass{amsart}
\usepackage{amssymb}
\usepackage{pdflscape}
\usepackage{amscd}
\usepackage{hyperref}
\usepackage{fancyhdr}

\newtheorem{theorem}{Theorem}[section]
\newtheorem{corollary}[theorem]{Corollary}
\newtheorem{lemma}[theorem]{Lemma}

\title{Local intersection cohomology of varieties of complexes}

\author{Xin Fang}
\address{Lehrstuhl f\"ur Algebra und Darstellungstheorie, RWTH Aachen, Pontdriesch 10-16, 52062 Aachen, Germany}
\email{xinfang.math@gmail.com}
\author{Markus Reineke}
\address{Ruhr-Universit\"at Bochum, Faculty of Mathematics, Universit\"atsstra{\ss}e 150, 44780 Bochum, Germany}
\email{markus.reineke@rub.de}

\begin{document}
\begin{abstract} We compute the local intersection cohomology of the irreducible components of varieties of complexes, by using Lusztig's geometric approach to quantum groups and explicit constructions of elements of Lusztig's canonical bases.\end{abstract}
\maketitle

\section{Introduction} The Buchsbaum-Eisenbud varieties of complexes \cite{BE}, parametrizing degree one differentials on finite-dimensional graded vector spaces, form an interesting and well-studied class of singular affine varieties. Their irreducible components are proven to be normal, Cohen-Macaulay, with rational singularities \cite{K,DL,DCS,LM} and Frobenius-split \cite{MT}, their defining equations, singular loci and their divisor class group are known \cite{DCS,L,Y}, and open embeddings into Schubert varieties are constructed \cite{MS}.\\[1ex]
In the present work, we complement this with quantitative topological information on the irreducible components, by describing in Theorem \ref{main} the Poincar\'e polynomials of their local intersection cohomology (that is, the stalks of cohomology sheaves of their intersection cohomology complexes).\\[1ex]
To achieve this, we use the geometric realization of canonical bases of quantized enveloping algebras \cite{LuCan,LuBook}, which identifies the relevant Poincaré polynomials with the base change coefficients between canonical and PBW type bases. The elements of the former corresponding to the irreducible components of varieties of complexes can be constructed using ideas of \cite{Xi}.\\[1ex]
We note that a similar strategy, obtaining new geometric information from calculations in quantized enveloping algebras, was employed in earlier work \cite{FR} in relation to families of degenerate flag varieties. We also note that the irreducible components of varieties of complexes do not admit small resolutions in general, making it more surprising that all local intersection cohomology can be computed explicitly.\\[1ex]
We review the definition of the variety of complexes in Section \ref{s1} and give a parametrization of their irreducible components. In Section \ref{s2}, we review the geometric realization of canonical bases of quantized enveloping algebras via Hall algebras. The relevant elements of the canonical basis are constructed in Section \ref{s3}, and the geometric consequences are worked out in Section \ref{s4}.

\vskip 10pt
\noindent
\textbf{Acknowledgments:}  The work of X.F. is funded by the Deutsche
Forschungsgemeinschaft: “Symbolic Tools in Mathematics and their Application” (TRR 195, project-ID 286237555).

\section{Varieties of complexes}\label{s1}

Let $k$ be an algebraically closed field (its characteristic will play a role later when considering $\ell$-adic, respectively singular, cohomology). Let $V_1,\ldots,V_n$ be finite-dimensional $k$-vector spaces of dimensions $d_1,\ldots,d_n$, respectively, and denote ${\bf d}=(d_1,\ldots,d_n)$.  Define ${\rm Com}(V_*)={\rm Com}({\bf d})$ as the affine variety
$$\{(f_1,\ldots,f_{n-1})\, |\, f_{i+1}\circ  f_i=0,\, i=1,\ldots,n-1\}\subset\bigoplus_{i=1}^{n-1}{\rm Hom}_k(V_i,V_{i+1}).$$

A point $f_*$ in ${\rm Com}(V_*)$ thus defines a complex of $k$-vector spaces
$$V_1\stackrel{f_1}{\rightarrow}V_2\stackrel{f_2}{\rightarrow}\ldots\stackrel{f_{n-1}}{\rightarrow}V_n.$$

The reductive algebraic group
$$G={\rm GL}(V_1)\times\ldots\times {\rm GL}(V_n)$$
acts on ${\rm Com}(V_*)$ via
$$(g_1,\ldots,g_n)\cdot(f_1,\ldots,f_{n-1})=(g_2f_1g_1^{-1},\ldots,g_nf_{n-1}g_{n-1}^{-1})$$
with finitely many orbits $\mathcal{O}({\bf r},{\bf h})$ for tuples of non-negative integers
$${\bf r}=(r_1,\ldots,r_{n-1}),\, {\bf h}=(h_1,\ldots,h_n)$$
such that
$$d_i=r_{i-1}+r_i+h_i,\, i=1,\ldots,n$$
(where we formally set $r_0=0=r_{n}$). Namely, the orbit $\mathcal{O}({\bf r},{\bf h})$ consists of all $f_*$ such that ${\rm rk}(f_i)=r_i$ for all $i=1,\ldots,n-1$ and (consequently)
$$\dim {\rm Ker}(f_i)-\dim{\rm Im}(f_{i-1})=h_i,\; i=1,\ldots,n$$
(where we formally set $f_0=0=f_n$).  In other words, the orbits are the loci of complexes with the same ranks and Betti numbers.\\[1ex]
We can also describe closure of these orbits, see \cite[§11]{MS}. Namely, we have $\mathcal{O}({\bf r}',{\bf h}')\subset \overline{\mathcal{O}({\bf r},{\bf h})}$ if and only if ${\bf r}'\leq{\bf r}$ componentwise. An orbit $\mathcal{O}({\bf r},{\bf h})$ is thus open if it is maximal with respect to this ordering, equivalently, if $h_ih_{i+1}=0$ for $i=1,\ldots,n-1$. We denote $$\Omega=\Omega({\bf r},{\bf h})=\{i=1,\ldots,n\, |\, h_i\not=0\}.$$
This set is thus sparse in the sense that $i\not\in \Omega$ or $i+1\not\in \Omega$ for $i=1,\ldots,n-1$. This proves:

\begin{lemma} The irreducible components of ${\rm Com}({\bf d})$ are the orbit closures $\overline{\mathcal{O}({\bf r},{\bf h})}$ for tuples $({\bf r},{\bf h})$ such that $\Omega({\bf r},{\bf h})$ is sparse.
\end{lemma}

\section{Quantized enveloping algebras}\label{s2}

\subsection{Quiver approach to canonical bases}

Let $C=(a_{ij})_{i,j\in I}$ be a finite type symmetric Cartan matrix. We consider the positive part $\mathcal{U}^+$ of the corresponding quantized enveloping algebra, which is the $\mathbb{Q}(v)$-algebra with Chevalley generators $E_i$ for $i\in I$ and defining relations
$$E_iE_j=E_jE_i\mbox{ if }a_{ij}=0,$$
$$E_i^2E_{j}-(v+v^{-1})E_iE_{j}E_i+E_{j}E_i^2=0\mbox{ if }a_{ij}=-1.$$
For $n\in\mathbb{N}$, we define quantum numbers and quantum factorials $$[n]=\frac{v^n-v^{-n}}{v-v^{-1}},\; [n]!=[1]\cdot[2]\cdot\ldots\cdot[n]$$ in $\mathbb{Z}[v,v^{-1}]$, as well as quantum binomial coefficients
$$\left[{a\atop n}\right]=\prod_{i=1}^n\frac{v^{a+1-i}-v^{-a-1+i}}{v^i-v^{-i}}\in\mathbb{Z}[v,v^{-1}]$$ for $a\in\mathbb{Z}$. Note that \cite[1.3.1.(a)]{LuBook}
$$\left[{a\atop n}\right]=(-1)^a\left[{{-a+n-1}\atop n}\right].$$
For $a\in\mathbb{N}$ and $q=v^2$, we also consider a normalized form
$$\left({a\atop n}\right)_{q}=v^{n(a-n)}\left[{a\atop n}\right]\in\mathbb{Z}[q].$$

For $x\in\mathcal{U}^+$, we define the divided power $x^{(n)}=\frac{1}{[n]!}x^n\in\mathcal{U}^+$.\\[1eX]
Let $Q$ be a quiver orienting the Dynkin diagram associated to $C$. Thus the set of vertices of $Q$ equals $I$, and there is an arrow between the vertices $i$ and $j$ if and only if $a_{ij}=-1$.  Define the Euler form of $Q$ as the non-symmetric bilinear form on $\mathbb{Z}I$ given by
$$\langle{\bf d},{\bf e}\rangle=\sum_{i\in I}d_ie_i-\sum_{\alpha:i\rightarrow j}d_ie_j;$$
its symmetrization is the quadratic form associated to $C$. The elements ${\bf d}\in \mathbb{N}I$ such that $\langle{\bf d},{\bf d}\rangle=1$ form the positive part $\Phi^+$ of the root system associated to $C$.\\[1ex]
We consider finite-dimensional representations of $Q$ over an arbitrary field $K$. By Gabriel's theorem, the isomorphism classes of indecomposable representations of $Q$ are in bijection to $\Phi^+$, and thus the isomorphism classes of arbitrary representations are in bijection to $\mathbb{N}$-valued functions on $\Phi^+$. For such a function $\alpha$, we denote by $M(\alpha,K)$ the corresponding $K$-representation of $Q$.\\[1ex]
For three such functions $\alpha,\beta,\gamma$, there exists a polynomial (called Hall polynomial) $F_{\alpha,\beta}^\gamma(q)\in\mathbb{Z}[q]$ such that, for any finite field $K$, the evaluation $F_{\alpha,\beta}^\gamma(|K|)$ equals the number of subrepresentations $U\subset M(\gamma,K)$ which are isomorphic to $M(\beta,K)$, with quotient $M(\gamma,K)/U$ isomorphic to $M(\alpha,K)$. By slight abuse of notation, we also denote $F_{\alpha,\beta}^\gamma(q)$ by $F_{M,N}^X(q)$ if $M=M(\alpha,K)$, $N=M(\beta,K)$, $X=M(\gamma,K)$ (again for an arbitrary field $K$).\\[1ex]
Similarly, for every function $\alpha$, there exists a polynomial $a_\alpha(q)\in\mathbb{Z}[q]$ such that $a_\alpha(|K|)$ equals the cardinality of ${\rm Aut}_Q(M(\alpha,K))$ for every finite field $K$. Again by abuse of notation we write $a_M(q)=a_\alpha(q)$ if $M=M(\alpha,K)$ for an arbitrary field $K$.\\[1ex]
For two representations $M,N$ of $Q$ over $K$, we abbreviate $\dim_K{\rm Hom}_K(M,N)$ by $[M,N]$.
We define the (generic, twisted) Hall algebra $H(Q)$ of $Q$ as the $\mathbb{Q}(v)$-vector space with basis elements $E_{[M]}$ indexed by the isomorphism classes of finite-dimensional representations of $Q$, and with multiplication defined by
$$E_{[M]}\cdot E_{[N]}=\sum_{[X]}v^{[M,M]+[N,N]+\langle{\rm\bf dim}\, M,{\rm\bf dim}\, N\rangle-[X,X]}F_{M,N}^X(v^2)E_{[X]}.$$
The algebras $\mathcal{U}^+$ and $H(Q)$ can be identified by mapping $E_i$ to $E_{[S_i]}$ for $i\in I$, where $S_i$ denotes the one-dimensional representation of $Q$ concentrated at the vertex $i$. The elements $E_{[M]}$ then form a basis of PBW type $B_Q=B_{\bf i}$ in the sense of \cite{LuCan}, for a reduced expression ${\bf i}$ of the longest element of the Weyl group of the Cartan matrix $C$ which is adapted to $Q$.\\[1ex]
Let $\overline{\cdot}$ be the unique $\mathbb{Q}$-linear involution on $\mathcal{U}^+$ fixing all $E_i$ and interchanging $v$ and $v^{-1}$. The canonical basis $\mathcal{B}$ of $\mathcal{U}^+$ is the unique $\mathbb{Z}[v^{-1}]$-basis of the $\mathbb{Z}[v^{-1}]$-lattice generated by $B_Q$ which is elementwise fixed by $\overline{\cdot}$. More concretely, it consists of elements $\mathcal{E}_{[M]}$ such that $\overline{\mathcal{E}_{[M]}}=\mathcal{E}_{[M]}$ and, in the expansion $$\mathcal{E}_{[M]}=\sum_{[N]}\zeta^M_NE_{[N]},$$ we have $\zeta^M_M=1$ and $\zeta^M_N\in v^{-1}\mathbb{Z}[v^{-1}]$ for all $N\not\simeq M$.\\[1ex]
Let now $k$ be an algebraic closure of a finite field $\mathbb{F}_q$. Fixing a dimension vector ${\bf d}\in\mathbb{N}I$ for $Q$, as well as $k$-vector spaces $V_i$ of dimension $d_i$ for $i\in I$, respectively, we consider the space of representations
$$R_{\bf d}(Q)=\bigoplus_{\alpha:i\rightarrow j}{\rm Hom}_k(V_i,V_j),$$ on which the group $G_{\bf d}=\prod_{i\in I}{\rm GL}(V_i)$ acts via base change $$(g_i)_i\cdot(f_\alpha)_\alpha=(g_jf_\alpha g_i^{-1})_{\alpha:i\rightarrow j}.$$

The orbits for this action are naturally in bijection to the isomorphism classes of representations of $Q$ of dimension vector ${\bf d}$; we denote the orbit corresponding to an isomorphism class $[M]$ by $\mathcal{O}_{[M]}$.\\[1ex]
Let $\ell$ be a prime invertible in $k$, and let ${\rm\bf IC}_{[M]}$ be the $\ell$-adic intersection cohomology complex on the closure $\overline{\mathcal{O}_{[M]}}$, extended trivially to a complex of constructible $\ell$-adic sheaves on $R_{\bf d}(Q)$. It is known \cite[Section 10]{LuCan} that ${\rm\bf IC}_{[M]}$ has strong purity properties. Namely, the stalk $\mathcal{H}^i_N({\rm\bf IC}_{[M]})$ at a point $N\in R_{\bf d}(Q)$ of the $i$-th cohomology sheaf of the complex ${\rm\bf IC}_{[M]}$ vanishes for $i$ odd, and the Frobenius map induced from $\mathbb{F}_q$ acts on it with all eigenvalues equal to $q^{i/2}$ in case $i$ is even.\\[1ex]
The main result on which our approach to the description of intersection cohomology is based is the fact \cite[Corollary 10.7]{LuCan} that the canonical basis and the local intersection cohomology  determine each other in the sense that
\begin{equation}\label{zetaic}v^{[N,N]-[M,M]}\zeta^M_N=\sum_i\dim\mathcal{H}^i_N({\rm\bf IC}_{[M]})v^i.\end{equation}

The same identity then holds over the base field $k=\mathbb{C}$ for singular intersection cohomology.







\subsection{Type $A$ quiver and complexes}

We specialize the previous constructions to the case of a Cartan matrix of type $A_n$ and the quiver 
$$1\rightarrow 2\rightarrow \ldots\rightarrow n.$$

The Euler form on $\mathbb{Z}I\simeq\mathbb{Z}^n$ is then given by
$$\langle{\bf d},{\bf e}\rangle=\sum_{i=1}^nd_ie_i-\sum_{i=1}^{n-1}d_ie_{i+1}.$$

The roots are $\alpha_{i,j}={\bf e}_i+\ldots+{\bf e}_j\in\mathbb{Z}I$ for $1\leq i\leq j\leq n$, corresponding to the indecomposable representations $U_{i,j}$ of $Q$ such that ${\rm\bf dim}\, U_{i,j}=\alpha_{i,j}$. Consequently, every representation $M$ of $Q$ is of the form
$$M\simeq \bigoplus_{1\leq i\leq j\leq n}U_{i,j}^{\oplus m_{i,j}}$$
for non-negative integers $m_{i,j}$ (encoding a function on positive roots as above).\\[1ex]
For a fixed dimension vector ${\bf d}\in\mathbb{N}I$ and vector spaces $V_i$ of dimension $d_i$ for $i\in I$, respectively, we have
$$R_{\bf d}(Q)=\bigoplus_{i=1}^{n-1}\mathrm{Hom}(V_i,V_{i+1}),$$
thus the variety of complexes $\mathrm{Com}(V_*)$ is naturally identified as a closed subvariety of $R_{\bf d}(Q)$. Under this identification, a representation $M$ of $Q$ belongs to ${\rm Com}(V_*)$ if and only if it is a direct sum of one- and two-dimensional indecomposables, thus 
$$M=M({\bf r},{\bf h})=\bigoplus_{i=1}^n S_i^{\oplus h_i}\oplus\bigoplus_{i=1}^{n-1}U_{i,i+1}^{\oplus r_i}.$$
for tuples ${\bf r}$ and ${\bf h}$ as before.\\[1ex]
The PBW type basis $B_Q$ of $\mathcal{U}^+$ equals the basis $B_{\bf i}$ for the  reduced expression $${\bf i}=(n,\, n-1,n,\, n-2,n-1,n,\, \ldots,\, 1,\ldots,n).$$
In $\mathcal{U}^+$, using the root elements $$E_{i,i+1}=E_iE_{i+1}-v^{-1}E_{i+1}E_i,$$
 we have
 $$E_{[M({\bf r},{\bf h})]}=E({\bf r},{\bf h})=E_n^{(h_n)}E_{n-1,n}^{(r_{n-1})}E_{n-1}^{(h_{n-1})}E_{n-2,n-1}^{(r_{n-2})}E_{n-2}^{(h_{n-2})}\ldots E_2^{(h_2)}E_{1,2}^{(r_1)}E_1^{(h_1)}.$$
 
Our aim is to construct the elements of the canonical basis of $\mathcal{U}^+$ corresponding to the elements $E({\bf r},{\bf h})$ with $\Omega({\bf r},{\bf h})$ sparse. More precisely, for $\mathbf{a}\in\mathbb{R}^n$ and $\mathbf{b}\in\mathbb{R}^{n-1}$, define ${\bf a}\dotplus{\bf b}$ as the tuple 
$$(a_1+b_1,a_2+b_1+b_2,\ldots,a_{n-1}+b_{n-2}+b_{n-1},a_n+b_{n-1}).$$
Then we want to construct elements $\mathcal{E}({\bf r},{\bf h})\in\mathcal{U}^+$ such that $\overline{\mathcal{E}({\bf r},{\bf h})}=\mathcal{E}({\bf r},{\bf h})$ and
$$\mathcal{E}({\bf r},{\bf h})=\sum_{{\bf k}\leq{\bf r}}\zeta_{\bf k}\cdot E({\bf r}-{\bf k},{\bf h}\dotplus{\bf k})$$
for elements $\zeta_{\bf k}$ such that $\zeta_{\bf  0}=1$ and $\zeta_{\bf k}\in v^{-1}\mathbb{Z}[v^{-1}]$ for all ${\bf k}\not={\bf 0}$, since the relevant intersection cohomology complex is supported on the closure $\overline{\mathcal{O}({\bf r},{\bf h})}$, which is the union of the $\mathcal{O}({\bf r}-{\bf k},{\bf h}\dotplus{\bf k})$.

\section{Construction of elements of canonical bases}\label{s3}

We consider certain elements of $\mathcal{U}^+$ and combine them in order to construct elements of the canonical basis. To make our constructions uniform for $i=1,\ldots,n$, we formally set $E_0=E_{0,1}=E_{n+1}=E_{n,n+1}=1$ in $\mathcal{U}^+$ (and $r_0=r_n=1$ as before).\\[1ex]
We continue to work with sequences ${\bf r},{\bf h}$ as above. For $1\leq i\leq j\leq n$, we define the element
$$B_{i,j}=B_{i,j}({\bf r},{\bf h})=E_{j-1}^{(r_{j-1})}E_j^{(r_{j-1})}\ldots E_{i+1}^{(r_{i+1})}E_{i+2}^{(r_{i+1})} E_i^{(r_i)}E_{i+1}^{(r_{i})}$$
(note that $B_{i,i}=1$). For $i=1,\ldots,n$, we define
 the element
$$C_i=C_i({\bf r},{\bf h})=\sum_{u_i=0}^{r_i}\left[{-h_i\atop u_i}\right]E_{i}^{(r_i-u_i)}E_{i+1}^{(r_i)}E_{i-1}^{(r_{i-1})}E_i^{(h_i+r_{i-1}+u_i)}.$$


Now assume that $\Omega=\Omega({\bf r},{\bf h})=\{i_1,\ldots,i_s\}$ with $i_1<\ldots<i_s$ is sparse. We then have in $\mathcal{U}^+$ a well-defined element 
$$\mathcal{E}_\Omega({\bf r},{\bf h})=B_{i_s+1,n}\cdot C_{i_s}\cdot B_{i_{s-1}+1,i_s-1}\cdot C_{i_{s-1}}\cdot\ldots\cdot C_{i_2}\cdot B_{i_1+1,i_2-1}\cdot C_{i_1}\cdot B_{1,i_1-1}.$$

\begin{theorem} Under the previous assumptions, $\mathcal{E}_\Omega({\bf r},{\bf h})$ belongs to $\mathcal{B}$.
\end{theorem}

We start with the identity \cite[Lemma 42.1.2.]{LuBook}
\begin{equation}\label{id1}E_i^{(m)}E_{i+1}^{(n)}=\sum_{k=0}^{m}v^{-(n-m+k)k}E_{i+1}^{(n-m+k)}E_{i,i+1}^{(m-k)}E_i^{(k)}\end{equation}
for $m\leq n$. Applying this repeatedly to the factors of $B_{i,j}$, we can rewrite $B_{i,j}$ as
$$\sum_{(k_p\leq r_p)_{i\leq p<j}}v^{-\sum_{i\leq p<j}k_p^2}E_j^{(k_{j-1})}E_{j-1,j}^{(r_{j-1}-k_{j-1})}E_{j-1}^{(k_{j-1})}\cdot E_{j-1}^{(k_{j-2})}E_{j-2,j-1}^{(r_{j-2}-k_{j-2})}E_{j-2}^{(k_{j-2})}\cdot\ldots\cdot$$
$$\cdot E_{i+2}^{(k_{i+1})}E_{i+1,i+2}^{(r_{i+1}-k_{i+1})}E_{i+1}^{(k_{i+1})}\cdot E_{i+1}^{(k_{i})}E_{i,i+1}^{(r_{i}-k_{i})}E_{i}^{(k_{i})}=$$
$$=\sum_{(k_p\leq r_p)_{i\leq p<j}}v^{-\sum_{i\leq p<j}k_p^2}\prod_{i<p<j}\left[{{k_{p-1}+k_p}\atop k_p}\right]E_j^{(k_{j-1})}E_{j-1,j}^{(r_{j-1}-k_{j-1})}\cdot$$
$$\cdot E_{j-1}^{(k_{j-2}+k_{j-1})}E_{j-2,j-1}^{(r_{j-2}-k_{j-2})}\cdot\ldots E_{i+1,i+2}^{(r_{i+1}-k_{i+1})}E_{i+1}^{(k_i+k_{i+1})}E_{i,i+1}^{(r_{i}-k_{i})}E_{i}^{(k_{i})}.$$
Applying identity (\ref{id1}) twice to the factors of $C_i$, we can rewrite $C_i$ as
$$\sum_{u_i=0}^{r_i}\left[{-h_i\atop u_i}\right]\left(\sum_{l_i=0}^{r_i-u_i}v^{-(u_i+l_i)l_i}E_{i+1}^{(u_i+l_i)}E_{i,i+1}^{(r_i-u_i-l_i)}E_i^{(l_i)}\right)\cdot$$
$$\cdot\left(\sum_{k_{i-1}=0}^{r_{i-1}}v^{-(h_i+u_i+k_{i-1})k_{i-1}}E_i^{(h_i+u_i+k_{i-1})}E_{i-1,i}^{(r_{i-1}-k_{i-1})}E_{i-1}^{(k_{i-1})}\right).$$
Multiplying the divided powers of $E_i$ and denoting $k_i=u_i+l_i$, this can be written as
\begin{equation}\label{t1}\sum_{k_{i-1}=0}^{r_{i-1}}\sum_{k_i=0}^{r_i}\left(\sum_{u_i}v^{-k_i(k_i-u_i)-(h_i+u_i+k_{i-1})k_{i-1}}\left[{-h_i\atop u_i}\right]\left[{h_i+k_i+k_{i-1}\atop k_i-u_i}\right]\right)\cdot\end{equation}
$$\cdot E_{i+1}^{(k_i)}E_{i,i+1}^{(r_i-k_i)} E_i^{(h_i+k_{i-1}+k_i)}E_{i-1,i}^{(r_{i-1}-k_{i-1})}E_{i-1}^{(k_{i-1})}.$$

We now use the following identity \cite[3.1.(b)]{Xi}, which holds for all $c_1,c_2,c_3\in\mathbb{Z}$ and $t\geq 0$:
$$\sum_{t=t_1+t_{23}}v^{(c_2-c_3)t_1-c_1t_{23}}\left[{c_1\atop t_1}\right]\left[{c_2+c_3\atop t_{23}}\right]=\sum_{t=t_{12}+t_3}v^{(c_2-c_1)t_3-c_3t_{12}}\left[{c_1+c_2\atop t_{12}}\right]\left[{c_3\atop t_3}\right].$$


To apply this identity, we rewrite the $v$-exponent in (\ref{t1}) as follows:
$$-k_i(k_i-u_i)-(h_i+u_i+k_{i-1})k_{i-1}=$$
$$=-h_i(k_{i-1}+k_i)-k_{i-1}^2-k_i^2+(h_i+k_i-k_{i-1})u_i-(-h_i)(k_i-u_i).$$
This allows us to apply the above identity with the parameters
$$c_1=-h_i,\, c_2=h_i+k_i,\, c_3=k_{i-1},\,   t=k_i,\, t_1=u_i,\, t_3=t_i$$
to obtain the following expression for $C_i$:
$$\sum_{k_{i-1}=0}^{r_{i-1}}\sum_{k_i=0}^{r_i}v^{-h_i(k_{i-1}+k_i)-k_{i-1}^2-k_i^2-k_{i-1}k_i}\left(\sum_{t_i}v^{(2h_i+k_{i-1}+k_i)t_i}\left[{k_i\atop k_i-t_i}\right]\left[{k_{i-1}\atop t_i}\right]\right)\cdot$$
$$\cdot E_{i+1}^{(k_i)}E_{i,i+1}^{(r_i-k_i)} E_i^{(h_i+k_{i-1}+k_i)}E_{i-1,i}^{(r_{i-1}-k_{i-1})}E_{i-1}^{(k_{i-1})}.$$

We now use these expansions of the elements $B_{i,j}$ and $C_i$ to expand all terms in the product

$$B_{i_s+1,n}\cdot C_{i_s}\cdot B_{i_{s-1}+1,i_s-1}\cdot C_{i_{s-1}}\cdot\ldots\cdot C_{i_2}\cdot B_{i_1+1,i_2-1}\cdot C_{i_1}\cdot B_{1,i_1-1}.$$

We see that this results in a sum over indices $k_i=0,\ldots,r_i$ for all $i=1,\ldots,n-1$, as well as indices $t_i=0,\ldots,\min(k_{i-1},k_i)$ whenever $i\in \Omega$. The terms in this expansion are products of the form
$$E_n^{(k_{n-1})}\ldots E_{i_s+1}^{(k_{i_s+1})}\cdot E_{i_s+1}^{(k_{i_s})}\ldots E_{i_s-1}^{(k_{i_s-1})}\cdot E_{i_s-1}^{(k_{i_s-2})}\ldots E_{i_{s-1}+1}^{(k_{i_{s-1}-1})}\cdot E_{i_{s-1}+1}^{(k_{i_{s-1}})}\ldots E_{i_{s-1}-1}^{(k_{i_{s-1}+1})}\cdot$$
$$\cdot\ldots\cdot E_{i_2+1}^{(k_{i_2})}\ldots E_{i_2-1}^{(k_{i_2-1})}\cdot E_{i_2-1}^{(k_{i_2-2})}\ldots E_{i_1+1}^{(k_{i_1+1})}\cdot E_{i_1+1}^{(k_{i_1})}\ldots E_{i_1-1}^{(k_{i_1-1})}\cdot E_{i_1-1}^{(k_{i_1-2})}\ldots E_1^{(k_1)}.$$

Multiplying divided powers of the $E_i$ thus produces additional factors
$$\left[{k_{i-1}+k_i\atop k_i}\right]$$
for all $i$ such that $i\pm 1\in \Omega$. We thus have the following expansion for $\mathcal{E}_\Omega({\bf r},{\bf h})$:
$$\sum_{(k_i\leq r_i)_i} \sum_{(t_i\leq\min(k_{i-1},k_i))_{i\in \Omega}} v^{-\sum_ik_i^2-\sum_{i\in \Omega}k_{i-1}k_i-\sum_{i\in \Omega}h_i(k_{i-1}+k_i)+\sum_{i\in \Omega}(2h_i+k_{i-1}+k_i)t_i}\cdot$$
$$\cdot\prod_{i\in \Omega}\left(\left[{k_i\atop k_i-t_i}\right]\left[{k_{i-1}\atop t_i}\right]\right)\cdot\prod_{i\not\in \Omega}\left[{k_{i-1}+k_i\atop k_i}\right]\cdot E({\bf r}-{\bf k},{\bf h}\dotplus{\bf k}).$$
We now estimate the degree of the coefficient to $E({\bf r}-{\bf k},{\bf h}\dotplus{\bf k})$. This degree equals
$$-\sum_ik_i^2-\sum_{i\in \Omega}k_{i-1}k_i-\sum_{i\in \Omega}h_i(k_{i-1}+k_i)+\sum_{i\in \Omega}(2h_i+k_{i-1}+k_i)t_i+$$
$$+\sum_{i\in \Omega}(k_i-t_i)t_i+\sum_{i\in \Omega}(k_{i-1}-t_i)t_i+\sum_{i\not\in \Omega}k_{i-1}k_i.$$
For each $i\in \Omega$, we rewrite
\begin{eqnarray*}
& &-k_{i-1}k_i-h_i(k_{i-1}+k_i)+(2h_i+k_{i-1}+k_i)t_i+(k_i-t_i)t_i+(k_{i-1}-t_i)t_i\\
&=& k_{i-1}k_i-h_i(k_{i-1}-t_i+k_i-t_i)-2(k_{i-1}-t_i)(k_i-t_i).
\end{eqnarray*}
This allows us to rewrite the degree as
$$-(\sum_ik_i^2-\sum_ik_{i-1}k_i)-\sum_{i\in \Omega}h_i(k_{i-1}-t_i+k_i-t_i)-2\sum_{i\in \Omega}(k_{i-1}-t_i)(k_i-t_i).$$
Every summand being nonnegative, we see that the degree is nonpositive, and equal to zero only if ${\bf k}={\bf 0}$. \\[1ex]
This finishes the proof of the theorem, since the element $\mathcal{E}_\Omega({\bf r},{\bf h})$ is obviously bar-invariant.

\begin{corollary}\label{zeta} The coefficient of $E({\bf r}-{\bf k},{\bf h}\dotplus{\bf k})$ in $\mathcal{E}_\Omega({\bf r},{\bf h})$ is given by
$$ \sum_{(t_i\leq\min(k_{i-1},k_i))_{i\in \Omega}} v^{-\sum_ik_i^2-\sum_{i\in \Omega}k_{i-1}k_i-\sum_{i\in \Omega}h_i(k_{i-1}+k_i)+\sum_{i\in \Omega}(2h_i+k_{i-1}+k_i)t_i}\cdot$$
$$\cdot\prod_{i\in \Omega}\left(\left[{k_i\atop k_i-t_i}\right]\left[{k_{i-1}\atop t_i}\right]\right)\cdot\prod_{i\not\in \Omega}\left[{k_{i-1}+k_i\atop k_i}\right].$$
\end{corollary}

\section{Application to local intersection cohomology}\label{s4}

Using the explicit expansion of the elements $\mathcal{E}_\Omega({\bf r},{\bf h})$ developed in the previous section, we can now easily deduce our main result.

\begin{theorem}\label{main} Let ${\bf r},{\bf h}$ be as before with $\Omega=\Omega({\bf r},{\bf h})$ sparse. Then the closure $\overline{\mathcal{O}({\bf r},{\bf h})}$ is an irreducible component of ${\rm Com}({\bf d})$. A stalk $\mathcal{H}^*_{f_*}({\rm\bf IC}(\overline{\mathcal{O}({\bf r},{\bf h})}))$ over a point $f_*\in\mathcal{O}({\bf r}-{\bf k},{\bf h}\dotplus{\bf k})$ has Poincar\'e polynomial
\begin{eqnarray*}
& &\sum_i\dim \mathcal{H}^{2i}_{f_*}({\rm\bf IC}(\overline{\mathcal{O}({\bf r},{\bf h})}))q^i\\
&=&\sum_{(t_i\leq\min(k_{i-1},k_i))_{i\in\Omega}}q^{\sum_{i\in\Omega}(h_i+t_i)t_i}\cdot\prod_{i\in\Omega}\left({k_i\atop k_i-t_i}\right)_q\left({k_{i-1}\atop t_i}\right)_q\cdot\prod_{i\not\in\Omega}\left({k_{i-1}+k_i\atop k_i}\right)_q.
\end{eqnarray*}
\end{theorem}

\proof Let $M=M({\bf r},{\bf h})$ and $N=M({\bf r}-{\bf k},{\bf h}\dotplus{\bf k})$. By the identity (\ref{zetaic}), the desired Poincar\'e polynomial can be expressed as
$$v^{[N,N]-[M,M]}\zeta^M_N,$$
where $\zeta^M_N$ is given by the formula in the previous corollary. Using the representation theory of the equioriented type $A$ quiver, we easily find: 
$$[M,M]=\sum_{i=1}^n(h_i^2+h_ir_i+r_i^2+h_{i}r_{i-1}+r_{i-1}r_{i}).$$

A lengthy cancellation then leads to 
$$[N,N]-[M,M]=\sum_{i=1}^{n}(h_ik_i+h_{i}k_{i-1}+k_i^2+k_{i-1}k_{i})$$
(where again we set $k_0=0=k_n$).\\[1ex]
Using furthermore the definition of the normalized quantum binomial coefficients, we see that Corollary \ref{zeta} implies the following formula for the Poincar\'e polynomial:

$$ \sum_{(t_i\leq\min(k_{i-1},k_i))_{i\in \Omega}} v^{*}\cdot\prod_{i\in \Omega}\left({k_i\atop k_i-t_i}\right)_{v^2}\left({k_{i-1}\atop t_i}\right)_{v^2}\cdot\prod_{i\not\in \Omega}\left({k_{i-1}+k_i\atop k_i}\right)_{v^2},$$
where the $v$-exponent $*$ equals
$$-\sum_{i=1}^nk_i^2-\sum_{i\in \Omega}k_{i-1}k_i-\sum_{i\in\Omega}h_i(k_{i-1}+k_i)+\sum_{i\in\Omega}(2h_i+k_{i-1}+k_i)t_i-$$
$$-\sum_{i\in \Omega}(k_i-t_i)t_i-\sum_{i\in \Omega}(k_{i-1}-t_i)t_i-\sum_{i\not\in\Omega}k_{i-1}k_i+$$
$$+\sum_{i=1}^nh_ik_i+\sum_{i=1}^nh_ik_{i-1}+\sum_{i=1}^nk_i^2+\sum_{i=1}^nk_{i-1}k_{i}=$$
$$=2\sum_{i\in \Omega}(h_i+t_i)t_i.$$
Replacing $v^2$ by $q$ yields the claimed formula.

\bibliography{ICSheafVC}{}
\bibliographystyle{plain}

\end{document}